\documentclass{article}%
\usepackage{amssymb}
\usepackage{amsmath}
\usepackage{amsfonts}
\usepackage{graphicx}%
\providecommand{\U}[1]{\protect\rule{.1in}{.1in}}
\newtheorem{theorem}{Theorem}

\newtheorem{example}[theorem]{Example}

\newtheorem{proposition}[theorem]{Proposition}

\newenvironment{proof}[1][Proof]{\noindent\textbf{#1.} }{\ \rule{0.5em}{0.5em}}
\begin{document}

\title{\textbf{SOME CHARACTERIZATIONS OF EULER SPIRALS IN }$E_{1}^{3}$}
\author{\textbf{Yusuf YAYLI}\thanks{Ankara University, Faculty of Science, Department
of Mathematics, Ankara, TURKEY}\textbf{, Semra SARACOGLU}\thanks{Siirt
University, Faculty of Science and Arts, Department of Mathematics, Siirt,
TURKEY}}
\maketitle

\begin{abstract}
In this study, some characterizations of Euler spirals in $E_{1}^{3}$ have
been presented by using their main property that their curvatures are linear.
Moreover, discussing some properties of Bertrand curves and helices, the
relationship between these special curves in $E_{1}^{3}$ have been
investigated with different theorems and examples. The approach we used in
this paper is useful in understanding the role of Euler spirals in $E_{1}^{3}$
in differential geometry.

\textbf{AMS Subj. Class.:} 53A04, 53A05, 53B30 .


\end{abstract}

\section{INTRODUCTION}

In three dimensional Euchlidean space $E^{3},$ Euler spirals are well-known as
the curves whose curvatures evolves linearly along the curve. It is also
called Clothoid or Cornu spiral whose curvature is equal to its arclength.

The equations of Euler spirals were written by Bernoulli first, in 1694. He
didn't compute these curves numerically. In 1744, Euler rediscovered the
curve's equations, described their properties, and derived a series expansion
to the curve's integrals. Later, in 1781, he also computed the spiral's end
points. The curves were re-discovered in 1890 for the third time by Talbot,
who used them to design railway tracks [1].

A new type of Euler spirals in $E^{2}$ and in $E^{3}$ are given in [1] with
their properties. They prove that their curve satisfies properties that
characterize fair and appealing curves and reduces to the 2D Euler spiral in
the planar case. Furthermore, they require that their curve conforms with the
definition of a 2D Euler spiral. Smiliarly, these curves are presented in [6]
as the ratio of two rational linear functions and have been defined in $E^{3}$
as generalized Euler spirals with some various characterizations. On the other
hand, linear relation between principal curvatures of spacelike surfaces in
Minkowski space is studied in [3].

In this paper, we present the timelike and spacelike Euler spirals in
Minkowski space $E_{1}^{3}.$ At first, we give the basic concepts and theorems
about the study then we deal with these spirals whose curvatures and torsion
are linear. Here, we seek that if any timelike Euler spirals in $E_{1}^{3}$ is
regular or not. Next, we investigate in which conditions the timelike Euler
spirals can be Bertrand curve. Additionally, by using the definition of
logarithmic spiral having a linear radius of curvature and a radius of torsion
from [1,2,6], we obtain the spacelike logarithmic spiral.

We believe that this study gives us a link and relation between the classical
differential surface theory and Euler spirals in Minkowski space $E_{1}^{3}.$

\section{PRELIMINARIES}

Now, we recall the basic concepts and important theorems on classical
differential geometry about the study, then we use them in the next sections
to give our approach. References [1,2,5,7], contain these concepts.

Let%

\[
\alpha:I\rightarrow E^{3}%
\]%
\[
\text{\ \ \ \ \ \ \ }s\mapsto\alpha(s)
\]
be non-null curve and $\left\{  T,N,B\right\}  $ frame of $\alpha.$ $T,N,B$
are the unit tangent, principal normal and binormal vectors respectively. Let
$\kappa$ and $\tau$ be the curvatures of the curve $\alpha.$

We consider a regular curve $\alpha$ parametrized by the length-arc. We call%
\[
\overrightarrow{T}(s)=\frac{d\alpha}{ds}%
\]
the tangent vector at $s.$ In particular, $\left\langle T(s),T^{\prime
}(s)\right\rangle =0.$ We will assume that $T^{\prime}(s)\neq0.$ In this study
we investigate the curve $\alpha$ in two different cases timelike and spacelike.

We suppose that $\alpha$ is a timelike curve. Then $T^{\prime}(s)\neq0$ is a
spacelike vector independent with $T(s).$ We define the curvature of $\alpha$
at $s$ as $\kappa(s)=\left\vert T^{\prime}(s)\right\vert .$The normal vector
$N(s)$ is defined by
\[
N(s)=\frac{T^{\prime}(s)}{\kappa(s)}=\frac{\alpha^{\prime\prime}%
(s)}{\left\vert \alpha^{\prime\prime}(s)\right\vert }.
\]
Moreover $\kappa(s)=\left\langle T^{\prime}(s),N(s)\right\rangle .$ We call
the binormal vector $B(s)$ as
\[
B(s)=T(s)\times N(s)
\]
The vector $B(s)$ is unitary and spacelike. For each $s,$ $\left\{
T,N,B\right\}  $ is an orthonormal base of $E_{1}^{3}$ which is called the
Frenet trihedron of $\alpha.$ We define the torsion of $\alpha.$ We define the
torsion of $\alpha$ at $s$ as%
\[
\tau(s)=\left\langle N^{\prime}(s),B(s)\right\rangle \text{ }[5].
\]
By differentiation each one of the vector functions of the Frenet trihedron
and putting in relation with the same Frenet base, we obtain the Frenet
equations, namely,%
\[
\left[
\begin{array}
[c]{c}%
T^{\prime}\\
N^{\prime}\\
B^{\prime}%
\end{array}
\right]  =\left[
\begin{array}
[c]{ccc}%
0 & \kappa & 0\\
\kappa & 0 & \tau\\
0 & -\tau & 0
\end{array}
\right]  \left[
\begin{array}
[c]{c}%
T\\
N\\
B
\end{array}
\right]
\]

On the other hand, Let $\alpha$ be a spacelike curve. These are three
possibilities depending on the casual character of $T^{\prime}(s)$

\textbf{1. }The vector $T^{\prime}(s)$ is spacelike. Again, we write
$\kappa(s)=\left\vert T^{\prime}(s)\right\vert ,$ $N(s)=T^{\prime}%
(s)/\kappa(s)$ and $B(s)=T(s)\times N(s).$ The vectors $N$ and $B$ are called
the normal vector and the binormal vector respectively. The curvature of
$\alpha$ is defined by $\kappa.$ The equations are
\[
\left[
\begin{array}
[c]{c}%
T^{\prime}\\
N^{\prime}\\
B^{\prime}%
\end{array}
\right]  =\left[
\begin{array}
[c]{ccc}%
0 & \kappa & 0\\
-\kappa & 0 & \tau\\
0 & \tau & 0
\end{array}
\right]  \left[
\begin{array}
[c]{c}%
T\\
N\\
B
\end{array}
\right]
\]
The torsion of $\alpha$ is defined by $\tau=-\left\langle N^{\prime
},B\right\rangle $.

\textbf{2.} The vector $T^{\prime}(s)$ is timelike. The normal vector is
$N=T^{\prime}/\kappa,$ where $\kappa(s)=\sqrt{-\left\langle T^{\prime
}(s),T^{\prime}(s)\right\rangle }$ is the curvature of $\alpha.$ The binormal
vector is $B(s)=T(s)\times N(s),$ which is spacelike vector. Now, the Frenet
equations
\[
\left[
\begin{array}
[c]{c}%
T^{\prime}\\
N^{\prime}\\
B^{\prime}%
\end{array}
\right]  =\left[
\begin{array}
[c]{ccc}%
0 & \kappa & 0\\
\kappa & 0 & \tau\\
0 & \tau & 0
\end{array}
\right]  \left[
\begin{array}
[c]{c}%
T\\
N\\
B
\end{array}
\right]
\]
The torsion of $\alpha$ is $\tau=-\left\langle N^{\prime},B\right\rangle $.

\textbf{3.} The vector $T^{\prime}(s)$ is lightlike for any $s$ (recall that
$T^{\prime}(s)\neq0$ and it is not proportional to $T(s)$). We define the
normal vector as $N(s)=T^{\prime}(s),$ which is independent linear with
$T(s).$ Let $B$ be the unique lightlike vector such that $\left\langle
N,B\right\rangle =1$ and orthogonal to $T.$ The vector $B(s)$ is called the
binormal vector of $\alpha$ at $s.$ The Frenet equations are%

\[
\left[
\begin{array}
[c]{c}%
T^{\prime}\\
N^{\prime}\\
B^{\prime}%
\end{array}
\right]  =\left[
\begin{array}
[c]{ccc}%
0 & 1 & 0\\
0 & \tau & 0\\
-1 & 0 & -\tau
\end{array}
\right]  \left[
\begin{array}
[c]{c}%
T\\
N\\
B
\end{array}
\right]
\]
The function $\tau$ is called the torsion of $\alpha.$ There is not a
definition of the curvature of $\alpha$ $[5].$

It is well-known that a planar curve of $E_{1}^{3}$ is included in a affine
plane and also this plane is a vector plane. In [5], planar curves with
constant curvature is studied. If the curve $\alpha$ is a planar curve in
$E^{2}$ parametrized by the length-arc and $\nu$ is a fixed unitary direction,
we call $\theta(s)$ the angle that makes $T(s)$ and $\nu,$ that is
\[
\cos(\theta(s))=\left\langle T(s),\nu\right\rangle \text{ }[5].
\]
It can be proved that the curve $\alpha$ is $\left\vert \theta^{\prime
}(s)\right\vert .$

Harary and Tall define the Euler spirals in $E^{3}$ the curve having both its
curvature and torsion evolve linearly along the curve$.$ Furthermore, they
require that their curve conforms with the definition of a $Cornu$ $spiral.$

Here, these curves whose curvatures and torsion evolve linearly are called
Euler spirals in $E_{1}^{3}$. Thus for some constants $a,b,c,d$ $\epsilon$ $%
\mathbb{R}
,$%
\begin{align}
\kappa(s)  &  =as+b\\
\tau(s)  &  =cs+d\nonumber
\end{align}
Next, we define logarithmic spiral having a linear radius of curvature and a
linear radius of torsion from [1,2]. They seek a spiral that has both a linear
radius of curvature and a linear radius of torsion in the arc-length
parametrization $s:$%
\begin{align*}
\kappa(s)  &  =\frac{1}{as+b}\\
\tau(s)  &  =\frac{1}{cs+d}%
\end{align*}
where $a,b,c$ and $d$ are constants.

Accordingly, the Euler spirals in $E_{1}^{2}$ and in $E_{1}^{3}$ that satisfy
equation above by a set of differential equations is the spacelike curve
$\alpha$ for which the following conditions hold:%
\begin{align*}
\frac{d\overrightarrow{T}(s)}{ds}  &  =\left(  \frac{1}{(as+b)}\right)
\overrightarrow{N}(s),\\
\frac{d\overrightarrow{N}(s)}{ds}  &  =\left(  \frac{1}{(as+b)}\right)
\overrightarrow{T}(s)+\left(  \frac{1}{(cs+d)}\right)  \overrightarrow{B}(s)\\
\frac{d\overrightarrow{B}(s)}{ds}  &  =-\left(  \frac{1}{(cs+d)}\right)
\overrightarrow{N}(s).
\end{align*}
In addition to these, we want to give our definition that Euler spirals in
$E_{1}^{3}$ whose ratio between its curvature and torsion evolve linearly is
called generalized Euler spirals in $E_{1}^{3}.$ Thus for some constants
$a,b,c,d$ $\epsilon$ $%
\mathbb{R}
,$%
\[
\dfrac{\kappa}{\tau}=\frac{as+b}{cs+d}%
\]

\begin{theorem}
Let $\alpha$ be a timelike curve parametrized by the length-arc and included
in a timelike plane. Let $\nu$ be a unit fixed vector of the plane pointing to
the future. Let $\phi(s)$ be the hyperbolic angle between $T(s)$ and $\nu.$
Then $\kappa(s)=\left\vert \phi^{\prime}(s)\right\vert $ $[5].$
\end{theorem}

\begin{theorem}
Let the curve $\alpha$ be a timelike planar Euler spirals in $E_{1}^{3}.$ Then
$\alpha$ is a Bertrand curve and there are two constants $A$ and $B;$ the
curvature $\kappa$ and the torsion $\tau$ such that
\[
A\kappa+B\tau=1\text{ \ \ }[5].
\]

\end{theorem}

\section{PLANAR EULER SPIRALS IN $E_{1}^{3}$}

In this section, we will discuss the main properties of timelike and spacelike
Euler spirals. At that time, the relationship between Euler spirals, Bertrand
curves, regular curves and helices are given with some theorems and cases.

\subsection{The Timelike Planar Euler Spirals in $E_{1}^{3}$}

Let $\phi:I\rightarrow%
\mathbb{R}
$ be the function. We assume that the plane is timelike and the plane $P$ be
timelike, that is $P=\left\langle E_{2},E_{3}\right\rangle $ and then let
parametrized curve by the arc-length be the curve $\beta.$ Thus, this curve
can be given as
\[
\beta(s)=y(s)E_{2}+z(s)E_{3}%
\]
with
\[
y^{\prime}(s)^{2}-z^{\prime}(s)^{2}=-1.
\]
From that,%
\[
\beta(s)=\left(  \overset{s}{\underset{s_{0}}{\int}}\cosh\phi(t)dt,\overset
{s}{\underset{s_{0}}{\int}}\sinh\phi(t)dt\right)
\]
with the linear curavature%
\[
\phi(s)=\overset{s}{\underset{s_{0}}{\int}}\kappa(u)du.
\]
This approach shows that a plane curve with any given smooth function as its
signed curvature can be found. But simple curvature can lead to these curves.

Thus, the function $\phi\in%
\mathbb{R}
$ such that
\begin{align*}
y^{\prime}(s)  &  =\sinh\phi(t)dt,\\
z^{\prime}(s)  &  =\cosh\phi(t)dt
\end{align*}
taking the derivative of this curve,%
\[
\beta^{\prime}(s)=\left(  \sinh\phi(s),\cosh\phi(s)\right)
\]
and giving the dot product
\[
\left\langle \beta^{\prime},\beta^{\prime}\right\rangle =\sinh\phi
(s)^{2}-\cosh\phi(s)^{2}=-1
\]
the timelike Euler spiral can be computed with the property%
\[
\kappa(s)=\left\vert \beta^{\prime\prime}(s)\right\vert =\left\vert
\phi^{\prime}(s)\right\vert =s.
\]
Therefore, taking the signed curvature be $\kappa(s)=s$ and $s_{0}=0,$ the
timelike Euler spiral $\beta$ can be obtained as
\[
\beta(s)=\left(  \overset{s}{\underset{0}{\int}}\cosh\frac{s^{2}}%
{2}ds,\overset{s}{\underset{0}{\int}}\sinh\frac{s^{2}}{2}ds\right)
\]

\begin{example}
\textbf{ }Let the signed curvature be $\kappa(s)=as+b.$ Taking $s_{0}=0,$ we
get
\begin{align*}
\phi^{\prime}(s)  &  =as+b\text{ and}\\
\phi(s)  &  =\frac{a}{2}s^{2}+bs.
\end{align*}
then the timelike planar curve $\beta$ can be given as%
\[
\beta(s)=\left(  \overset{s}{\underset{0}{\int}}\cosh(\frac{a}{2}%
s^{2}+bs)ds,\overset{s}{\underset{0}{\int}}\sinh(\frac{a}{2}s^{2}%
+bs)ds\right)  .
\]

\end{example}

\begin{example}
\textbf{ }Let the signed curvature be $\kappa(s)=\dfrac{1}{s}.$ Taking
$s_{0}=0,$ we get the timelike planar curve as%
\[
\beta(s)=\left(  \overset{s}{\underset{0}{\int}}\cosh(\ln s)ds,\overset
{s}{\underset{0}{\int}}\sinh(\ln s)ds\right)
\]

\end{example}

\begin{theorem}
Any timelike planar Euler spirals in $E_{1}^{3}$ is regular.
\end{theorem}

\begin{proof}
Assume that the curve $\beta$ is timelike planar Euler spiral and is written
by%
\[
\beta(s)=\left(  \overset{s}{\underset{s_{0}}{\int}}\cosh\phi(s)ds,\overset
{s}{\underset{s_{0}}{\int}}\sinh\phi(s)ds\right)
\]
where
\[
\beta^{\prime}(s)=\left(  x^{\prime}(s),y^{\prime}(s),z^{\prime}(s)\right)
=(0,y^{\prime}(s),z^{\prime}(s))
\]
and%
\[
\left\langle \beta^{\prime}(s),\beta^{\prime}(s)\right\rangle =y^{\prime
}(s)^{2}-z^{\prime}(s)^{2}\langle0,
\]
In particular $\overset{s}{\underset{s_{0}}{%
{\displaystyle\int}
}}\sinh\phi(s)\neq0$ that is, $\beta$ is a regular curve.
\end{proof}

\subsection{The Spacelike Planar Euler Spirals in $E_{1}^{3}$}

Let $\phi:I\rightarrow%
\mathbb{R}
$ be the function. We suppose that the plane is timelike and is given by
$\left\{  x=0\right\}  .$ We are going to find the spacelike planar curve
$\alpha;$ at this time, the curve $\alpha$ is written as
\begin{align*}
\alpha(s)  &  =x(s)E_{1}+z(s)E_{3}\text{ and }\\
\alpha(s)  &  =(x(s),0,z(s))
\end{align*}
with%
\[
x^{\prime}(s)^{2}-z^{\prime}(s)^{2}=1.
\]
Thus, the spacelike planar Euler spiral is given by
\[
\alpha(s)=\left(  \overset{s}{\underset{s_{0}}{\int}}\sinh\phi(t)dt,\overset
{s}{\underset{s_{0}}{\int}}\cosh\phi(t)dt\right)
\]
where
\[
\phi(s)=\overset{s}{\underset{s_{0}}{\int}}\kappa(u)du
\]
with the curvature $\kappa.$ Here, as it is known the curvature $\kappa$ is
linear and the function $\phi\in%
\mathbb{R}
$ such that
\begin{align*}
x^{\prime}(s)  &  =\cosh\phi(s)\text{ and}\\
z^{\prime}(s)  &  =\sinh\phi(s).
\end{align*}
Then, the tangent vector $\alpha$ is
\begin{align*}
\alpha^{\prime}(s)  &  =\left(  x^{\prime}(s),z^{\prime}(s)\right)  \text{
and}\\
\alpha^{\prime}(s)  &  =(\cosh\phi(s),\sinh\phi(s)).
\end{align*}
It is clear that
\[
\left\langle \alpha^{\prime}(s),\alpha^{\prime}(s)\right\rangle =\cosh
\phi(s)^{2}-\sinh\phi(s)^{2}=1.
\]
Let the signed curvature be $\kappa(s)=s.$ Thus, the spacelike planar Euler
spiral in $E_{1}^{3}$ \ can be found as
\[
\alpha(s)=\left(  \overset{s}{\underset{0}{\int}}\sinh\frac{s^{2}}%
{2}ds,\overset{s}{\underset{0}{\int}}\cosh\frac{s^{2}}{2}ds\right)  .
\]
\textbf{ }

\begin{example}
If we take
\[
\kappa(s)=as+b
\]
then
\[
\phi^{\prime}(s)=as+b
\]
and%
\[
\phi(s)=a\frac{s^{2}}{2}ds+bs
\]
So,%
\[
\alpha(s)=\left(  \overset{s}{\underset{0}{\int}}\sinh(\frac{a}{2}%
s^{2}+bs)ds,\overset{s}{\underset{0}{\int}}\cosh(\frac{a}{2}s^{2}%
+bs)ds\right)
\]

\end{example}

\begin{example}
If the curvature $\kappa$ is constant as $\kappa=a,$ then $\phi(s)=as+b.$
Thus, the spacelike planar curve is obtained as%
\[
\alpha(s)=\left(  \overset{s}{\underset{0}{\int}}\sinh(as+b)ds,\overset
{s}{\underset{0}{\int}}\cosh(as+b)ds\right)  .
\]

\end{example}

\begin{example}
If the curvature $\kappa$ is given as
\[
\kappa=\frac{1}{s},
\]
then
\[
\phi(s)=\ln s.
\]
Therefore, the spacelike planar curve is
\[
\alpha(s)=\left(  \overset{s}{\underset{0}{\int}}\sinh(\ln s)ds,\overset
{s}{\underset{0}{\int}}\cosh(\ln s)ds\right)  .
\]

\end{example}

\section{EULER SPIRALS IN $E_{1}^{3}$}

In this section, we study some characterizations of Euler spirals in
$E_{1}^{3}$ by giving some theorems and definitions from $[4,5,6]$.

\begin{proposition}
If the curvature $\tau$ is zero then $\kappa=as+b$ and also the curve is
planar cornu spiral in $E_{1}^{3}$.

\begin{proof}
If $\tau=0$ and the curvature is linear, then the ratio%
\[
\frac{\tau}{\kappa}=0
\]
Therefore, we see that the curve is planar Euler spiral in $E_{1}^{3}$.
\end{proof}
\end{proposition}

\begin{proposition}
If the curvatures are
\begin{align*}
\tau &  =as+b\\
\kappa &  =c
\end{align*}
then the Euler spirals are rectifying curves in $E_{1}^{3}$.

\begin{proof}
If we take the ratio%
\[
\frac{\tau}{\kappa}=\frac{as+b}{c}%
\]
where $\lambda_{1}$ and $\lambda_{2}$ , with $\lambda_{1}\neq0$ are constants,
then [4]
\[
\frac{\tau}{\kappa}=\lambda_{1}s+\lambda_{2}.
\]
It shows us that the Euler spirals are rectifying curves in $E_{1}^{3}$. It
can be easily seen from [3] that rectifying curves have very simple
characterization in terms of the ratio $\dfrac{\tau}{\kappa}.$
\end{proof}
\end{proposition}

\begin{proposition}
Euler spirals in $E_{1}^{3}$are Bertrand curves.

\begin{proof}
From the definition of Euler spiral in $E_{1}^{3}$ and the equations above, we
can take%
\begin{align*}
\tau(s)  &  =c_{1}s+c_{2}\\
\kappa(s)  &  =d_{1}s+d_{2},\text{ with }c_{1}\neq0\text{ and }d_{1}\neq0.
\end{align*}
Here,%
\[
s=\frac{1}{c_{1}}(\tau-c_{2})
\]
and then,%
\begin{align*}
\kappa &  =\frac{d_{1}}{c_{1}}(\tau-c_{2})+d_{2}\\
c_{1}\kappa &  =d_{1}(\tau-c_{2})+c_{1}d_{2}\\
c_{1}\kappa-d_{1}\tau &  =c_{3}\\
\frac{c_{1}}{c_{3}}\kappa-\frac{d_{1}}{c_{3}}\tau &  =1
\end{align*}
Thus, we obtain%
\[
\lambda\kappa+\mu\tau=1
\]
from \textbf{Theorem.2 }that Euler spirals are Bertrand curves in $E_{1}^{3}$.
\end{proof}
\end{proposition}

\begin{theorem}
Let $M$ and $M_{r}$ be parallel surfaces in $E_{1}^{3}$ and also let the curve
$\alpha$ be a geodesic timelike Euler spiral on the surface $M$ such that the
curvatures
\begin{align*}
\kappa(s)  &  =c_{1}s+c_{2}\\
\tau(s)  &  =d_{1}s+d_{2},\text{ with }c_{1}\neq0\text{ and }d_{1}\neq0.
\end{align*}
In this case, the Bertrand pair of the curve $\alpha$ is on the surface
$M_{r}.$ Here,%
\[
r=\frac{c_{1}}{c_{1}d_{2}-d_{1}c_{2}}.
\]

\end{theorem}

\begin{proof}
From the proposition before, the Euler spiral is Bertrand curve. Let $\beta$
be the Bertrand pair of the curve of $\alpha$ as:%
\[
\beta(s)=\alpha(s)+\nu N(s).
\]
From the property that $\alpha$ is geodesic on $M$%
\[
N(s)=n(s).
\]
Therefore,%
\[
\beta(s)=\alpha(s)+\nu n(s).
\]
From the \textbf{Proposition 12,} we can give
\[
r=\frac{c_{1}}{c_{1}d_{2}-d_{1}c_{2}}.
\]
Thus, if we take $r=\nu$ then we have $\beta(s)\in M_{r}.$
\end{proof}

\begin{theorem}
Let $M$ be a surface in $E_{1}^{3}$ and $\alpha:I\rightarrow M$ be non-null
curve (timelike or spacelike). If the Darboux curve
\[
W(s)=\varepsilon\tau T+\kappa B
\]
is geodesic curve on the surface $M,$ then the curve $\alpha$ is Euler spiral
in $E_{1}^{3}.$ Here, if the curve $\alpha$ is timelike then $\varepsilon=-1$
and if the curve $\alpha$ is spacelike then $\varepsilon=1.$
\end{theorem}

\begin{proof}
Let $\alpha$ be timelike. Thus,%
\[
W(s)=-\tau T-\kappa B
\]
then we have
\begin{align}
W^{\prime}(s)  &  =-\tau^{\prime}T-\kappa^{\prime}B\\
W^{\prime\prime}(s)  &  =-\tau^{\prime\prime}T-\kappa^{\prime\prime}%
B-(\kappa\tau^{\prime}-\tau\kappa^{\prime})N
\end{align}
and also
\[
W^{\prime\prime}(s)=\lambda(s)n(s).
\]
Here, $n(s)$ is the unit normal vector field of the surface $M.$ Darboux curve
is geodesic on surface $M$, therefore we have
\[
W^{\prime\prime}(s)=\lambda(s)N(s).
\]
and then it can be easily given that $n=N.$ If we take
\[
\tau^{\prime\prime}=0,\text{ }\kappa^{\prime\prime}=0
\]
and also%
\begin{align*}
\kappa &  =as+b\\
\tau &  =cs+d
\end{align*}
then the curve $\alpha$ is generalized Euler spiral in $E_{1}^{3}.$
\end{proof}

On the other hand, let $\alpha$ be spacelike. Thus,%
\[
W(s)=\tau T-\kappa B.
\]
Here, if the vector $T^{\prime}(s)$ is spacelike or timelike then we have
\begin{align*}
W^{\prime}(s)  &  =\tau^{\prime}T-\kappa^{\prime}B\\
W^{\prime\prime}(s)  &  =\tau^{\prime\prime}T-\kappa^{\prime\prime}%
B+(\kappa\tau^{\prime}-\tau\kappa^{\prime})N
\end{align*}
and also it can be seen that
\[
W^{\prime\prime}(s)=\lambda(s)n(s).
\]
Smiliarly, $n(s)$ is the unit normal vector field of the surface $M.$ Darboux
curve is geodesic on surface $M,$ therefore we have
\[
W^{\prime\prime}(s)=\lambda(s)N(s).
\]
and because of that it can be easily given $n=N.$ If we take
\[
\tau^{\prime\prime}=0,\text{ }\kappa^{\prime\prime}=0
\]
and also%
\begin{align*}
\kappa &  =as+b\\
\tau &  =cs+d
\end{align*}
then the curve $\alpha$ is generalized Euler spiral in $E_{1}^{3}.$

\section{GENERALIZED EULER SPIRALS IN $E_{1}^{3}$}

In this section, we investigate generalized Euler spirals in $E_{1}^{3}$ by
using the definitions in above$.$

\begin{theorem}
In $E_{1}^{3},$ all logarithmic spirals are generalized Euler spirals.

\begin{proof}
As it is known that in all logarithmic spirals, the curvatures are linear as:%
\begin{align*}
\kappa(s)  &  =\frac{1}{as+b}\\
\tau(s)  &  =\frac{1}{cs+d}%
\end{align*}
In that case, it is clear that the ratio between the curvatures can be given
as:%
\[
\dfrac{\kappa}{\tau}=\frac{cs+d}{as+b}%
\]
Thus, it can be easily seen all logarithmic spirals are generalized euler spirals.

\begin{theorem}
Euler spirals are generalized Euler spirals in $E_{1}^{3}.$

\begin{proof}
It is clear from the property of curvature, torsion and the ratio that are
linear as:%
\begin{align*}
\kappa(s)  &  =as+b\\
\tau(s)  &  =cs+d
\end{align*}
and then%
\[
\dfrac{\kappa}{\tau}=\frac{as+b}{cs+d}%
\]
That shows us Euler spirals in $E_{1}^{3}$ are generalized Euler spirals in
$E_{1}^{3}$.
\end{proof}
\end{theorem}
\end{proof}
\end{theorem}

\begin{proposition}
All generalized Euler spirals in $E_{1}^{3}$ that have the property
\[
\dfrac{\tau}{\kappa}=d_{1}s+d_{2}%
\]
are rectifying curves.

\begin{proof}
If the curvatures $\kappa(s)$ and $\tau(s)$ are taken as
\begin{align*}
\kappa(s)  &  =c\\
\tau(s)  &  =d_{1}s+d_{2}\text{ with }d_{1}\neq0
\end{align*}
then, from [4]%
\[
\dfrac{\tau}{\kappa}=\frac{d_{1}s+d_{2}}{c}=\lambda_{1}s+\lambda_{2}\text{ }%
\]
where $\lambda_{1}$ and $\lambda_{2}$ are constants. This gives us that if the
curve $\alpha$ is generalized Euler spiral in $E_{1}^{3}$ then it is also in
rectifying plane.
\end{proof}
\end{proposition}

\textbf{Result.1 }General helices are generalized Euler spirals in $E_{1}^{3}$.

\begin{proof}
It \ can be seen from the property of curvatures that are linear and the ratio
is also constant as it is shown:%
\[
\dfrac{\tau}{\kappa}=\lambda
\]

\end{proof}

\begin{theorem}
Let
\begin{equation}
\alpha:I\rightarrow E_{1}^{3}%
\end{equation}%
\[
\text{\ \ \ \ \ \ \ }s\mapsto\alpha(s)
\]
be non-null curve (spacelike or timelike) and let $\kappa$ and $\tau$ be the
curvatures of the Frenet vectors of the curve $\alpha.$ For $a,b,c,d,\lambda
\in%
\mathbb{R}
,$ let take the curve $\beta$ as%
\begin{equation}
\beta(s)=\alpha(s)+(as+b)T+(cs+d)B+\lambda N.
\end{equation}
In this case, the curve $\alpha$ is generalized Euler spiral in $E_{1}^{3}$
which has the property%
\[
\frac{\kappa}{\tau}=\varepsilon\frac{cs+d}{as+b},\text{ \ \ \ \ \ \ \ }%
\varepsilon=\mp1
\]
if and only if the curves $\beta$ and $(T)$ are the involute-evolute pair.
Here, the curve $(T)$ is the tangent indicatrix of $\alpha$.

\begin{proof}
The tangent of the curve $\beta$ is%
\begin{equation}
\beta^{\prime}(s)=((1-\lambda)\kappa+a)T+(c+\lambda\tau)B+(\kappa
(as+b)-\tau(cs+d))N.
\end{equation}
The tangent of the curve $(T)$ is%
\[
\frac{dT}{ds_{T}}=\varepsilon N.
\]
Here, $s_{T}$ is the arc parameter of the curve $(T).$%
\begin{equation}
\left\langle \beta^{\prime},N\right\rangle =\kappa(as+b)-\tau(cs+d)
\end{equation}
If the curves $\beta$ and $(T)$ are the involute-evolute pair then
\[
\left\langle \beta^{\prime},N\right\rangle =0.
\]
From (7), it can be easily obtained
\[
\frac{\kappa}{\tau}=\frac{cs+d}{as+b}%
\]
This means that the curve $\alpha$ is generalized Euler spiral in $E_{1}^{3}.$

On the other hand, if the curve $\alpha$ is generalized Euler spiral in
$E_{1}^{3}$ which has the property%
\[
\frac{\kappa}{\tau}=\frac{cs+d}{as+b},
\]
then from (7)%
\[
\left\langle \beta^{\prime},N\right\rangle =0.
\]
This means that $\beta$ and $(T)$ are the involute-evolute pair in $E_{1}%
^{3}.$
\end{proof}
\end{theorem}

\textbf{Result.2 }From the hypothesis of the theorem above, the curve $\alpha$
is generalized Euler spiral in $E_{1}^{3}$ which has the property%
\[
\frac{\kappa}{\tau}=\frac{cs+d}{as+b}%
\]
if and only if $\beta$ and $(B)$ are the involute-evolute pair. Here, $(B)$ is
the binormal of the curve $\alpha$ in $E_{1}^{3}.$

\begin{proof}
The tangent of the curve $(B)$ is
\[
\frac{dB}{ds_{B}}=-N.
\]
It can be easily seen from [5] for $(B)$.

\begin{theorem}
Let the ruled surface $\Phi$ be%
\begin{align}
\Phi &  :I\times%
\mathbb{R}
\longrightarrow E_{1}^{3}\\
(s,v)  &  \rightarrow\Phi(s,v)=\alpha(s)+v[(as+b)T+(cs+d)B]\nonumber
\end{align}
and the curve $\alpha:I\rightarrow M$ be non-null curve. The ruled surface
$\Phi$ is developable if and only if the curve $\alpha$ is generalized Euler
spiral in $E_{1}^{3}$ which has the property%
\[
\frac{\kappa}{\tau}=\varepsilon\frac{cs+d}{as+b},\text{ \ \ \ \ \ \ \ }%
\varepsilon=\mp1.
\]
Here, $\alpha$ is the base curve and $T,$ $B$ are the tangent and binormal of
the curve $\alpha,$ respectively.
\end{theorem}
\end{proof}

\begin{proof}
For the directrix of the surface
\[
X(s)=(as+b)T+(cs+d)B,
\]
and also for
\[
X^{\prime}(s)=aT+[(as+b)\kappa\pm(cs+d)\tau]N+cB,
\]
we can easily give that%
\[
\det(T,X,X^{\prime})=\left\vert
\begin{array}
[c]{ccc}%
1 & 0 & 0\\
as+b & 0 & cs+d\\
a & (as+b)\kappa\pm\tau(cs+d) & c
\end{array}
\right\vert
\]
In this case, the ruled surface is developable if and only if $\det
(T,X,X^{\prime})=0$ then%
\[
(cs+d)(as+b)\kappa-(cs+d)\tau=0
\]
then for $cs+d\neq0$%
\[
(as+b)\kappa\pm(cs+d)\tau=0.
\]
Thus, the curve $\alpha$ is generalized Euler spiral in $E_{1}^{3}$ which has
\ the property%
\[
\frac{\tau}{\kappa}=\frac{as+b}{cs+d}.
\]

\end{proof}

\begin{theorem}
Let $\alpha:I\rightarrow M$ be non-null curve. If the curve
\[
U(s)=\frac{\varepsilon}{\kappa}T-\frac{1}{\tau}B
\]
is a geodesic curve then the curve $\alpha$ is a logarithmic spiral in
$E_{1}^{3}.$
\end{theorem}

\begin{proof}
Let $\alpha$ be timelike. Thus we have
\[
U(s)=-\frac{1}{\kappa}T-\frac{1}{\tau}B
\]
then
\begin{align}
U^{\prime}  &  =\left(  -\frac{1}{\kappa}\right)  ^{\prime}T-\left(  \frac
{1}{\tau}\right)  ^{\prime}B\\
U^{\prime\prime}  &  =\left(  -\frac{1}{\kappa}\right)  ^{\prime\prime
}T-\left(  \frac{1}{\tau}\right)  ^{\prime\prime}B-\left[  \left(  \frac
{1}{\kappa}\right)  ^{\prime}\kappa-\left(  \frac{1}{\tau}\right)  ^{\prime
}\tau\right]  N.
\end{align}
$U(s)$ is geodesic on surface $M,$ therefore%
\[
U^{\prime\prime}=\mu_{1}(s)n(s).
\]
Here, $n(s)$ is the unit normal vector field of the surface $M.$ And also
$n=N,$ then we have
\[
U^{\prime\prime}=\mu_{1}(s)N(s).
\]
Also, from (10) it can be easily seen that
\begin{align*}
\frac{1}{\kappa}  &  =as+b\\
\frac{1}{\tau}  &  =cs+d
\end{align*}
Thus, it is clear that the curve $\alpha$ is a logarithmic spiral in
$E_{1}^{3}$.

On the other hand, let $\alpha$ be spacelike. Thus we have
\[
U(s)=\frac{1}{\kappa}T-\frac{1}{\tau}B
\]
then
\begin{align}
U^{\prime}  &  =\left(  \frac{1}{\kappa}\right)  ^{\prime}T-\left(  \frac
{1}{\tau}\right)  ^{\prime}B\\
U^{\prime\prime}  &  =\left(  \frac{1}{\kappa}\right)  ^{\prime\prime
}T-\left(  \frac{1}{\tau}\right)  ^{\prime\prime}B+\left[  \left(  \frac
{1}{\kappa}\right)  ^{\prime}\kappa-\left(  \frac{1}{\tau}\right)  ^{\prime
}\tau\right]  N.
\end{align}
$U(s)$ is geodesic on surface $M,$ therefore%
\[
U^{\prime\prime}=\mu_{2}(s)n(s).
\]
Here, $n(s)$ is the unit normal vector field of the surface $M.$ And also
$n=N,$ then we have
\[
U^{\prime\prime}=\mu_{2}(s)N(s).
\]
Also, from (12) it can be easily seen that
\begin{align*}
\frac{1}{\kappa}  &  =as+b\\
\frac{1}{\tau}  &  =cs+d
\end{align*}

\end{proof}

\section{CONCLUSIONS}

In this study, at the beginning, planar Euler spirals in $E_{1}^{3}$ have been
defined and then generalized Euler spirals in $E_{1}^{3}$ have been introduced
by using their important properties. Depending on these, some different
characterizations of Euler spirals in $E_{1}^{3}$ are expressed by giving
theorems, propositions with their results and proofs. At this time, it is
obtained that Euler spirals are generalized Euler spirals in $E_{1}^{3}.$
Additionally, we show that all logarithmic spirals are generalized Euler
spirals in $E_{1}^{3}.$ Moreover, many different approaches about generalized
Euler spirals in $E_{1}^{3}$ are presented in this paper.

We hope that this study will gain different interpretation to the other
studies in this field.

\end{document}